\documentclass[12pt,english]{article}
\usepackage{theorem}
\usepackage{a4wide}
\usepackage{graphics}
\usepackage{epsfig}
\usepackage{amssymb}
\usepackage{latexsym}
\usepackage{amsmath}

\newtheorem{thh}{Theorem}[section]
\newtheorem{df}[thh]{Definition}

\newtheorem{lem}[thh]{Lemma}

\newtheorem{prop}[thh]{Proposition}

\title{Cohomology of regular differential \\ forms for
affine curves}
\author{Philippe Bonnet}
\date{}

\newcommand{\dem}{{\em Proof: }}
\newcommand{\qed}{\begin{flushright} $\blacksquare$\end{flushright}}

\newcommand{\CC}{\mathbb C}

\newcommand{\CN}{\mathbb C ^n}
\newcommand{\AC}{{\cal{O}}_C}
\newcommand{\OS}{{\cal{O}}_S}
\newcommand{\CM}{{\cal{M}}}
\newcommand{\OX}{{\cal{O}}_X}

\newcommand{\ACN}{{\cal{O}}_{\widetilde{C}}}
\newcommand{\ACV}{{\cal{O}}_{C,V}}
\newcommand{\COCE}{\widehat{\Omega^1 _{C,V}}}
\newcommand{\OCE}{\Omega^1 _{C,V}}
\newcommand{\ACNV}{{\cal{O}}_{\widetilde{C},\widetilde{V}}}
\newcommand{\NC}{\widetilde{C}}
\newcommand{\NV}{\widetilde{V}}

\begin{document}
\maketitle

\begin{center} { \small
Mathematisches Institut, Universit\" at Basel\\
Rheinsprung 21, 4051 Basel, Switzerland.\\
E-mail:Philippe.bonnet@unibas.ch}
\end{center}

\begin{abstract}
Let $C$ be a complex affine reduced curve, and denote by $H^1(C)$ its first truncated cohomology group,
i.e. the quotient of all regular differential 1-forms by exact 1-forms. First we introduce
a nonnegative invariant $\mu'(C,x)$ that measures the complexity of the
singularity of $C$ at the point $x$, and we establish the following formula:
$$
dim \; H^1(C) = dim \; H_1(C) + \sum_{x \in C}\mu'(C,x)
$$
where $H_1(C)$ denotes the first singular homology group of $C$ with complex coefficients.
Second we consider a family of curves given by the fibres of a dominant morphism
$f:X\rightarrow \CC$, where $X$ is an irreducible complex affine surface. We
analyse the behaviour of the function $y\mapsto dim \; H^1(f^{-1}(y))$. More precisely
we show that it is constant on a Zariski open set, and that it is lower semi-continuous
in general.
\end{abstract}

\section{Introduction}

Let $C$ be a reduced complex affine curve that may be reducible or singular. For
any integer $k$, denote by $\Omega ^ k (C)$ the
space of regular differential $k$-forms (or K\"ahler forms) on $C$. The exterior
derivative $d$ is well-defined on $\Omega ^ k (C)$, and yields a complex:
$$
0\longrightarrow \CC \longrightarrow \Omega ^ 0 (C) \longrightarrow \Omega^1 (C)
\longrightarrow 0
$$
The first truncated De Rham cohomology group $H^1(C)$ is the
quotient $\Omega ^ 1 (C)/d\Omega ^ 0 (C)$.
If $C$ is smooth, then $C$ is a non-compact Riemann surface, for
which the De Rham cohomology
groups $H^k _{DR}(C)$ with complex coefficients are well-defined. Moreover
$H^1(C)$ coincides with the algebraic De Rham cohomology group of $C$ (see \cite{Ei})
and, by a theorem of Grothendieck (see \cite{Gro}), we have the isomorphism:
$$
H^1(C)\simeq H^1 _{DR}(C)
$$
So truncated De Rham cohomology is always defined
and coincides with standard De Rham cohomology if
$C$ is smooth. We would like to know to what extend this cohomology
reflects the topological properties of $C$, especially when $C$ has
singularities.

\begin{df}
Let $\widehat{\Omega ^k _{C,x}}$ be the space of
formal differential $k$-forms on the germ $(C,x)$.
The local De Rham cohomology group of $C$ at $x$
is the quotient:
$$
H^1 (C,x)= \widehat{\Omega ^1 _{C,x}}/d\widehat{\Omega ^ 0 _{C,x}}
$$
Its dimension $\mu'(C,x)$ is the local Betti number of $C$ at $x$.
\end{df}
This number characterizes the presence of singularities, in the sense that $\mu'(C,x)=0$
if and only if $x$ is a smooth point of $C$. Moreover it coincides
with the Milnor number (see \cite{Mi}) if $C$ is locally a complete intersection (see
\cite{B-G}).

Let $H_{1}(C)$ be the first singular homology group of $C$ with complex coefficients. By the results of Bloom and Herrera (see \cite{Bl-H}),
the integration of differential 1-forms along 1-cycles
is well-defined and provides us with a bilinear pairing $<,>$ on $H^1(C)\times H_1(C)$ given by:
$$<\omega,\gamma>= \int_{\gamma} \omega$$  
This induces the so-called De Rham morphism $\beta: H^1(C)\rightarrow H_1(C)^*, \; \omega \mapsto <\omega, .>$.By
Poincar\'e Duality and a theorem of Grothendieck (see \cite{Gro}), this map is an isomorphism when $C$ is smooth.
In the general case, we establish the following formula.

\begin{thh}\label{first} For any complex affine curve $C$, we have:
$dim \; H^1 (C)= dim\; H_1 (C) +
\sum_{x\in C} \mu'(C,x)$.
\end{thh}
The idea of the proof is the following. For any affine curve $C$, the morphism $\beta$ is onto (see \cite{Bl-H})
and this yields the exact sequence:
$$
0\longrightarrow \ker \beta \longrightarrow H^1(C) \longrightarrow H_1(C) ^* \longrightarrow 0
$$
For any point $x$ in $C$, every regular 1-form $\omega$ can be seen as a formal 1-form on the germ $(C,x)$.
Moreover every exact 1-form on $C$ is exact as a formal 1-form on $(C,x)$. We then have a natural morphism:
$$
i_x: H^1 (C) \longrightarrow H^1(C,x)
$$
We prove that the morphism $\alpha$:
$$
\alpha: \ker\;\beta \longrightarrow \oplus_{x\in C} H^1(C,x) \quad , \quad \omega\longmapsto
(i_x(\omega))_{x\in C}
$$
is an isomorphism, which gives the result by passing to the dimensions. \\

So local Betti numbers measure the default to Poincar\'e Duality
in the case of singular curves. Theorem \ref{first} implies in particular
that {\em a complex affine curve is isomorphic to a disjoint union of
copies of $\CC$ if and only if $H^1 (C)= 0$}. \\

Now we are going to study the behaviour of the function $h_1(y)=dim \; H^1(f^{-1}(y))$, where $X$
is a complex affine irreducible surface and $f: X\longrightarrow \CC$ is a dominant morphism.{\em The following
results still hold for any reduced surface $X$ (that is, any equidimensional reduced affine variety
of dimension 2) as soon as the morphism $f$ is dominant on every irreducible component of $X$}.
Recall that ${\cal{P}}$ holds for every generic point of $\CC$ if the set of points $y$ of $\CC$
where ${\cal{P}}(y)$ does not hold is finite. We have the following first result.

\begin{prop} \label{fini}
Let $X$ be a complex affine irreducible surface and $f: X \rightarrow\CC$ a dominant morphism. Then
there exists an integer $h_f\geq 0$ such that, for every generic point $y$ of $\CC$:
$$
dim \; H^1(f^{-1}(y))= h_f
$$
\end{prop}
For the proof, we introduce the truncated relative cohomology group $H^1(f)$ of $f$. We first show that it is finitely
generated after a suitable localisation. This is analogous to what happens for algebraic relative cohomology groups
(see \cite{Ha}, and
\cite{A-B} in the smooth case). If $h_f$ denotes the rank of $H^1(f)$ as a $\CC[f]$-module, we then show that it coincides with
the dimension of $H^1(f^{-1}(y))$ for generic $y$. 

\begin{thh} \label{second} 
Let $X$ be a complex affine surface that is locally a complete intersection and $f: X\rightarrow \CC$ be a dominant morphism.
If $f^{-1}(y)\cap Sing(f)$ is finite, then:
$$
dim \; H^1(f^{-1}(y)) \leq  h_f
$$
\end{thh}
In particular the function $h^1$ is lower semi-continuous at every point $y_0$ of $\CC$ such that $f^{-1}(y_0)\cap Sing(f)$ is finite, i.e:
$$
h^1(y_0)\leq {\underline{\lim}}_{y\to y_0} h^1(y)
$$
The previous results have analogous settings in terms of singular homology. Indeed if $X$ is equal to $\CC^2$ and
$f:\CC^2 \rightarrow \CC$ is a polynomial mapping, then there exists a non-empty
Zariski open set $U$ in $\CC$ such that $f:f^{-1}(U)\rightarrow U$ is a locally trivial topological fibration (see \cite{V}).
In particular all the fibres $f^{-1}(y)$, $y\in U$, are homeomorphic and there exists an integer $p$ such that:
$$
dim \; H_1(f^{-1}(y))= p
$$
for any $y$ in $U$. If $f$ has isolated singularities, then its local Betti number at any singular point coincides with
its Milnor number. If $f$ is moreover tame, then every fibre $f^{-1}(y)$ has the homotopy type of a bouquet of $\mu -\mu^y$
circles, where $\mu$ is the sum of all Milnor numbers and $\mu^y$ is the sum of the Milnor numbers of critical points lying on
$f^{-1}(y)$ (see \cite{Bro}). In particular this implies that:
$$
dim \; H_1(f^{-1}(y))= \mu -\mu^y
$$
By theorem \ref{first}, the dimension of $H^1(f^{-1}(y))$ is equal to $\mu$ if $f$ is tame, hence it is independent of $y$. If
$f$ is not tame, then the dimension of $H_1(f^{-1}(y))$ may be $<\mu -\mu_y$. This loss in the topology can be interpreted as
follows. The map $f$ possesses singularities at infinity where a certain number of cycles of the general fibre vanish (see \cite{S-T}).
In this case,
we see that the function $h^1(y)$ is lower semi-continuous, as asserted by theorem \ref{second}.

We end up this paper with an example of a mapping $f:X\rightarrow \CC$, where $X$ is not locally a complete
intersection. In this example, the dimension of $H^1(f^{-1}(y))$ increases for a special fibre,
which is contrary to what is predicted by theorem \ref{second}. This phenomenon is due to the presence
of a special singularity on $X$, which produces a local Betti number for $f^{-1}(0)$ where there should be
none.

\section{Properties of the normalisation}

Let $C$ be a complex affine curve and $\AC$ its ring of regular functions. Let
$\NC$ be its affine normalisation and $\Pi: \NC\rightarrow C$ the
normalisation morphism. Assume that $\NC$ is embedded in $\CN$, and let $\overline{B(0,R)}$ be the closed ball
of $\CN$ for the standard hermitian metric. For $R$ large enough, the intersection $S=\NC \cap \overline{B(0,R)}$
contains all the preimages of singular points of $C$ and is a deformation retract of $\NC$. We fix $R$
together with a triangulation $\widetilde{T}$ of $S$. Let $\widetilde{V}$ be its set of vertices. Since
$\Pi$ is finite, we may refine $\widetilde{T}$ so that the set $V=\Pi(\widetilde{V})$ contains all the singular
points of $C$, and so that $\widetilde{V}=\Pi^{-1}(V)$. By construction, the image
$T=\Pi(\widetilde{T})$ defines a triangulation of $\Pi(S)$. Since $\Pi$ is an isomorphism
from $\NC - \widetilde{V}$ to $C-V$, the set $\Pi(S)$ is also a deformation retract of $C$.
In particular, every 1-cycle of $C$ can be represented by a formal sum of edges of the triangulation
$T$. We denote by $\{\widetilde{\gamma_i}\}$
the set of edges of $\widetilde{T}$, and set $\gamma_i=\Pi(\widetilde{\gamma_i})$. {\em We consider
this triangulation fixed from now on}.
 
For any point $x$ in $C$, ${\cal{O}}_{C,x}$ stands for the ring of germs of regular functions at $x$. Denote
by $\ACV$ the ring of germs of regular functions at $V$, i.e. the direct sum:
$$
\ACV= \oplus_{x \in V} \; {\cal{O}}_{C,x}
$$
Let $I$ be the vanishing ideal of the set $V$ in $C$, and denote by $\widehat{\ACV}$ the $I$-adic
completion of $\ACV$. Note that we have the isomorphism:
$$
\widehat{\ACV}= \oplus_{x \in V} \; \widehat{{\cal{O}}_{C,x}}
$$
A {\em formal function on $(C,V)$} is an element of $\widehat{\ACV}$. In a similar way,
denote by $\Omega^1 _{C,x}$ the space of germs of regular 1-forms on $C$
at $x$, and by $\Omega^1 _{C,V}$ the finite sum:
$$
\Omega^1 _{C,V}= \oplus_{x \in V} \; \Omega^1 _{C,x}
$$
The $I$-adic completion $\COCE$ of $\OCE$ is the set of {\em formal 1-forms on $(C,V)$}. Note that
we have the isomorphism:
$$
\COCE= \oplus_{x \in V} \; \widehat{\Omega^1 _{C,x}}
$$
We can define the sets of {\em formal functions and formal 1-forms on} $(\NC,\widetilde{V})$ in exactly the same way. In this section,
we are going to
describe the relationships between the functions and 1-forms on $\NC$ and $C$.  

\subsection{Formal functions}

Let $\Pi^* : \AC \longrightarrow \ACN$ be the morphism induced by the normalisation map.
After localisation at $V$ and completion, we obtain the following injective map:
$$
\widehat{\Pi^* _V} : \widehat{\ACV} \longrightarrow \widehat{\ACNV}
$$
Since the germ $(\NC,x)$ is smooth for any point $x$ in $\NC$, every element $R$ of $\ACNV$
has a well-defined order $ord_x(R)$ at $x$, and thus it defines a divisor:
$$
div (R)=\sum_{x \in \NV} ord_{x}(R) x
$$
 
\begin{prop} \label{order1}
Let $\widetilde{R}$ be a formal function on $(\NC,\NV)$
that vanishes at every point of $\NV$. Then there exists a regular function
$S$ on $\NC$, vanishing at every point of $\NV$, and a formal function $R$
on $(C,V)$ such that $\widetilde{R}= S + \widehat{\Pi^* _V}(R)$.
\end{prop} 
In order to prove this proposition, we need the following lemma.
\begin{lem} \label{formal}
With the previous notations, there exists a divisor
$D$ on $(\NC,\NV)$ such that, for any formal function $\widetilde{R}$ on
$(\NC,\NV)$, we have: $div (\widetilde{R})\geq D \;\Rightarrow \; \widetilde{R}
\in \widehat{\Pi^* _V}(\widehat{\ACV})$.
\end{lem}
\dem Let $A$ be a conductor of the normalisation, i.e.
an element of $\AC$ that is not a zero-divisor and such that $\Pi^*(A)\ACN \subseteq \Pi^*(\AC)$.
After localisation at $V$ and completion, we obtain that:
$$
\widehat{\Pi^* _V} (A )\widehat{\ACNV}
\subseteq \widehat{\Pi^* _V}(\widehat{\ACV})
$$  
Set $D= div \; \widehat{\Pi^* _V}(A)$ and let $\widetilde{R}$ be a formal function
on $(\NC,\NV)$ such that $div (\widetilde{R})\geq D$. Then $\widetilde{R}$
is locally divisible by $\widehat{\Pi^* _V}(A)$, and the quotient $S=\widetilde{R}/\widehat{\Pi^* _V}(A)$
is a formal function on $(\NC,\NV)$. Therefore
$\widetilde{R} = \widehat{\Pi^* _V}(A)S$ belongs to $\widehat{\Pi^* _V}(\widehat{\ACV})$.
\qed
{\it Proof of Proposition \ref{order1}}: Let $\widetilde{R}$ be a formal function
on $(\NC,\NV)$. For any point $x$ in $\NV$, let $z_x$ be a uniformising parameter of
$\NC$ at $x$ defined on all of $\NC$. Then $\widetilde{R}$ has a Taylor expansion
$\sum_{k\geq 0} R_{k,x} z_x ^k$ at $x$. For any such $x$, we set:
$$
R_x= \sum_{k\leq n} R_{k,x} z_x ^k
$$
Let ${\cal{O}}_{\NC}$ be the ring of regular functions on $\NC$, and denote by $\widetilde{I}$ the ideal generated by $I$
in ${\cal{O}}_{\NC}$.  Since the radical of $\widetilde{I}$ is the vanishing ideal of $\widetilde{V}$, $\widehat{\ACNV}$ is the
$\widetilde{I}$-adic completion of ${\cal{O}}_{\NC}$. So there exists a regular function $S$ on $\NC$,
whose Taylor expansion of order $n$ at any point $x$ is equal to $R_x$.
For $n$ large enough, we have the inequality:
$$
div (\widetilde{R}-S) \geq D
$$
By lemma \ref{formal}, there exists a formal
function $R$ on $(C,V)$ such that $\widehat{\Pi^* _V}(R)=\widetilde{R} -S$.
\qed

\subsection{Formal 1-forms}

Let $\Pi^*:\Omega^1 (C)\rightarrow \Omega^1 (\NC)$ be the morphism induced by normalisation.
After localisation at $V$ and completion, we obtain the following morphism:
$$
\widehat{\Pi^* _V}: \widehat{\Omega^1 _{C,V}}\longrightarrow \widehat{\Omega^1 _{\NC,\NV}}
$$
In this subsection, we consider $\Omega^1 (\NC)$ as a $\AC$-module via the multiplication rule
$(P,\omega)\mapsto \Pi^*(P)\omega$. If $M$ is an $\AC$-module and $\CM$ is an ideal, denote by
$M_{\CM}$ its localisation with respect to $\CM$, and by $\widehat{M_{\CM}}$ its $\CM$-adic completion. 
We are going to prove the following proposition.

\begin{prop} \label{formm}
Let $\omega$ be a formal 1-form on the germ $(C,V)$. Then there exist a formal function
$R$ on $(C,V)$, a regular 1-form $\omega_0$ on $C$ and a regular function $S$ in $\ACN$,
vanishing at all points of $\NV$, such that
$\omega=dR + \omega_0$ and $\Pi^*(\omega_0)=dS$.
\end{prop}

\begin{lem} \label{com}
Let $R$ be a noetherian ring, and $L:M \rightarrow N$ a morphism of finite $R$-modules. Let
$\omega$ be an element of $N$ that belongs to $Im \; \widehat{L_{\cal{M}}}$
for any maximal ideal ${\cal{M}}$. Then $\omega$ belongs to $Im \; L$.
\end{lem}
\dem First we show that $\omega$ belongs to $Im \; L_{\cal{M}}$ for any maximal ideal
${\cal{M}}$. Let $\{e_1,..,e_k\}$ be a set of generators of $M$, i.e.
$M=R<e_1,..,e_k>$. After localisation and completion, we get the equalities:
$$
\widehat{M _{\cal{M}}}=\widehat{R_{\cal{M}}}<e_1,..,e_k> \quad \mbox{and} \quad
Im \; \widehat{L_{\cal{M}}}=\widehat{R_{\cal{M}}}<L(e_1),..,L(e_k)>
=\widehat{Im \; L_{\cal{M}}}
$$
Since $N$ has finite type, the ${\cal{M}}$-adic topology on $N$ is Hausdorff
and we find:
$$
Im \; L_{\cal{M}}= Im \; \widehat{L_{\cal{M}}} \cap N
$$
So $\omega$ belongs to $Im \; L_{\cal{M}}$, and for any maximal ideal ${\cal{M}}$, there
exists an element $P_{\cal{M}}$ of $R -{\cal{M}}$ such that $P_{\cal{M}} \omega$ belongs
to $Im \; L$. Let $I$ be the
ideal in $R$ generated by all the $P_{\cal{M}}$. We claim that $I=(1)$, so that
$\omega$ belongs to $Im \; L$. Indeed if $I$ were not equal to $(1)$, it would be
contained in a maximal ideal ${\cal{M}}_0$ by Zorn's Lemma. Since $I$ contains
$P_{{\cal{M}}_0}$, $P_{{\cal{M}}_0}$ would be contained in ${\cal{M}}_0$, hence
a contradiction.
\qed

\begin{lem} \label{surjec2}
Let $\widetilde{\omega}$ be an element of $\Omega^1(\NC)\cap Im \; \widehat{\Pi^* _{V}}$. Then
$\widetilde{\omega}$
belongs to $Im \; \Pi^*$. 
\end{lem}
\dem We set $M=\Omega^1(C)$, $N=\Omega^1(\NC)$ and $L=\Pi^*$. Let ${\cal{M}}$ be a maximal ideal
and $x$ the corresponding point in $C$. If $x$ belongs to $V$, then $\widetilde{\omega}$
belongs to $Im \; \widehat{L_{{\cal{M}}}}$ by assumption. If not, then $\widetilde{\omega}$ still
belongs to $Im \; \widehat{L_{{\cal{M}}}}$ because $x$ is a smooth point of $C$, and
then $\widehat{L_{{\cal{M}}}}$ is an isomorphism. By lemma \ref{com}, $\widetilde{\omega}$ belongs to
$Im \; \Pi^*$.
\qed

\begin{lem} \label{surjec3}
Under the previous assumptions, $dim \; ker \Pi^*$ is finite and the natural map $ker \Pi^*
\rightarrow ker \widehat{\Pi^* _{V}}$ is an isomorphism.
\end{lem}
\dem For any $x$ in $C$, denote by $\CM$ the vanishing ideal of $x$ and set $L=\Pi^*$. For any $x$ outside $V$, $\Pi$ is an isomorphism over an open
neighborhood of $x$. So the map $\widehat{L_{\CM}}$ is an isomorphism for all
$x$ outside $V$, and the support of $ker \Pi^*$ is contained in $V$. Since $V$ is a finite set and $ker \Pi^* $ is a
finite module, $ker \Pi^* $ is an artinian module and $dim \; ker \Pi^*<\infty$. So there exists an order $n$ such
that $I^{n}\; ker \Pi^*=0$, and $ker \Pi^*$ is complete for the $I$-adic topology. Since completion is an exact
functor, we have:
$$
ker \Pi^*\simeq \widehat{ker \Pi^*} \simeq ker \widehat{\Pi^* _{V}}
$$
\qed
{\it Proof of proposition \ref{formm}}: Let $\omega$ be a formal 1-form on the germ $(C, V)$.
Since the germ $(\NC,\NV)$ is a disjoint union of smooth curves, the 1-form
$\widehat{\Pi^* _{V}}(\omega)$ is exact on each of these curves. There exists a
formal function $\widetilde{R}$ on $(\NC,\NV)$ such that:
$$
\widehat{\Pi^* _{V}}(\omega)=d\widetilde{R}
$$
By proposition \ref{order1}, there exist a regular function $S$ on $\NC$, vanishing at all
points of $\NV$,
and a formal function $R$ on $(C,V)$ such that $\widetilde{R}=S +\widehat{\Pi^* _{V}}(R)$.
After derivation, this implies:
$$
\widehat{\Pi^* _{V}}(\omega-d\widetilde{R})=dS
$$
By lemma \ref{surjec2} applied to $\widetilde{\omega}=dS$, there exists a regular 1-form
$\omega_1$ on $C$ such that $\Pi^*(\omega_1)=dS$. This yields:
$$
\widehat{\Pi^* _{V}}(\omega-d\widetilde{R}-\omega_1)=0
$$
By lemma \ref{surjec3}, there exists a regular 1-form $\omega_2$ in $ker \Pi^*$ such that
$\omega-d\widetilde{R}-\omega_1=\omega_2$. Then the 1-form $\omega_0=\omega_1 + \omega_2$
is regular on $C$ and satisfies the following relations:
$$
\omega=d\widetilde{R}+ \omega_0 \quad \mbox{and} \quad \Pi^*(\omega_0)=dS
$$
\qed

\section{Proof of theorem \ref{first}}

Let $C$ be a complex reduced affine curve in $\CN$, and let $\beta: H^1(C)\rightarrow H_1(C)^*$ be the map
defined in the introduction. Since $\beta$ is onto, it induces the following complex:
$$
0\longrightarrow ker \; \beta  \longrightarrow H^1 (C)
\longrightarrow H_1(C)^*\longrightarrow 0
$$
Moreover the inclusion of regular 1-forms into formal 1-forms at $x$ induces
a morphism:
$$
\alpha: ker \; \beta \longrightarrow \oplus_{x \in C} H^1(C,x)
$$
Since $C$ carries a structure a $CW$-complex, the vector space $H_1(C)$ is finite dimensionnal, and the same holds for
every $H^1(C,x)$ (see \cite{B-G}). So for the proof of theorem \ref{first}, we only need to show that $\alpha$ is an
isomorphism, and the result will follow by passing to the dimensions.

\subsection{Injectivity of $\alpha$}

Without loss of generality, we may assume that the curve $C$ is connected. Let $\omega$ be an element of $ker \; \beta$.
Fix a point $x_0$ in $C$, and consider the map $R$ defined as follows. For any point $x$ in $C$, choose a path $\gamma$
going from $x_0$ to $x$, and set:
$$
R(x)= \int_{\gamma} \omega
$$
Since $\omega$ has null integral along any closed path in $C$, this number is
well-defined and independent of the path $\gamma$ chosen. Furthermore the
function $S=R\circ \Pi$ is holomorphic on $\NC$ because it defines
an integral of $\Pi^*(\omega)$ on $\NC$. By Grothendieck's Theorem,
$S$ is a regular function on $\NC$, and $S$ takes the value $R(x)$
on $\Pi^{-1}(x)$. 

Assume now that $\alpha(\omega)=0$. Then for any point $x$ of $C$, the class of $\omega$ in $H^1(C,x)$ is zero, and there exists a formal
function $R^x$ on the germ $(C,x)$ such that $\omega= dR^x$. Let $\CM$ be the vanishing ideal of $x$ and denote by $\widehat{L_{\CM}}$
the morphism induced by $\Pi^*$ after localisation at $\CM$ and completion. The formal function $S - \widehat{L_{\CM}}(R^x)$ on
$(\NC, \Pi^{-1}(x))$ is constant around every point of $\Pi^{-1}(x)$, because $S$ and
$\widehat{L_{\CM}}(R^x)$ are both integrals of $\Pi^*(\omega)$. Since $S$ and $\widehat{L_{\CM}}(R^x)$
are constant on $\Pi^{-1}(x)$, there exists a constant
$\lambda$ such that:
$$
S - \widehat{L_{\CM}}(R^x) =\lambda
$$
on $(\NC, \Pi^{-1}(x))$. Up to replacing $R^x$ by $R^x -\lambda$, we may assume that $\lambda=0$,
and so $S$ belongs to $Im \; \widehat{L_{\CM}}$ for any point $x$
in $C$. By applying lemma \ref{com} to the morphism $\Pi^* : \AC \rightarrow \ACN$
of finite $\AC$-modules, we get that $S$ belongs to $\AC$. Since $S= R^x$ for any
$x$ in $C$, we get by derivation:
$$
\omega= dS = dR^x \quad \mbox{in} \quad \widehat{\Omega^1 _{C,x}}
$$
Since $\Omega^1 _{C,x}$ is a finite ${\cal{O}}_{C,x}$-module, the ${\cal{M}}$-adic
topology is separated and $\omega= dS$ in $\Omega^1 _{C,x}$. By Bourbaki result (
Commutative Algebra Chap 1-7 Corollary 1, p. 88), $\omega= dS$ in $\Omega^1 (C)$
and the class of $\omega$ in $H^1(C)$ is zero.

\subsection{Surjectivity of $\alpha$}

By construction, the set $V$ contains all the singular points of $C$. Since $H^1(C,x)=0$ if $C$ is smooth at $x$, we have the isomorphism:
$$
\oplus_{x \in C} H^1(C,x)\simeq \oplus_{x \in V} H^1(C,x)
$$
So every element $\omega$ of this sum can be represented by a formal 1-form on $(C,V)$,
which we also denote by $\omega$. By lemma \ref{formm}, there exist a formal function $R$ on
$(C,V)$, a regular
1-form $\omega_0$ on $C$ and a regular function $S$ on $\NC$, vanishing at all points
of $\NV$, such that:
$$
\omega= dR + \omega_0 \quad \mbox{and} \quad \Pi^*(\omega_0)=dS 
$$
Let $\gamma$ be a 1-cycle in $C$. This cycle can be represented as a formal linear combination of the edges $\gamma_i$
of the triangulation $T$. Since $S$ vanishes at all vertices of $\widetilde{T}$, and these vertices are endpoints
of the $\widetilde{\gamma_i}$, we have:
$$
\int_{\gamma_i} \omega_0= \int_{\widetilde{\gamma_i}} \Pi^*(\omega_0)=
\int_{\widetilde{\gamma_i}} dS= S(\widetilde{\gamma_i}(1))-
S(\widetilde{\gamma_i}(0))=0
$$
By linearity, we get that $<\omega _0, \gamma>=0$ for any cycle $\gamma$ in $C$. So $\omega_0$ belongs to $ker \; \beta$ and represents
the same class as $\omega$ in $\oplus_{x \in V} H^1(C,x)$. Therefore $\alpha(\omega_0)=\omega$ and $\alpha$ is surjective.

\section{Relative cohomology}

Let $X$ be a complex irreducible affine surface, and $f: X\longrightarrow \CC$ a dominant morphism.
Denote by $\Omega^k(X)$ the space of regular
$k$-forms on $X$. The {\em first group
of truncated relative cohomology of $f$} is the quotient:
$$
H^1 (f)= \frac{\Omega^1(X)}{d\Omega^0(X) + \Omega^0(X)df}
$$
Note that $H^1(f)$ is a $\CC[f]$-module via the multiplication
$(P(f),\omega)\mapsto P(f)\omega$. In the case of analytic germs $f$,
relative cohomology groups have been extensively used to describe the topological
and cohomological properties of $f$; for more details, see for instance
\cite{Loo}. In the algebraic setting, the relative cohomology of polynomial
mappings has been intensively studied, especially via the use of the Gauss-Manin
connexion (see for instance \cite{A-B}).We are going to study some properties of
truncated relative cohomology and use them to prove proposition \ref{fini}.

\subsection{Finiteness of truncated relative cohomology}

In this subsection, we are going to establish that $H^1(f)$ is, after a suitable localisation,
a finite module. More precisely:

\begin{prop} \label{finite3}
Let $f:X\longrightarrow \CC$ be a dominant morphism, where $X$ is an irreducible surface. Then there
exists a non-zero polynomial $P$ of $\CC[t]$ such that $H^1(f)_{(P(f))}$ is a $\CC[f]_{(P(f))}$-module
of finite type.
\end{prop}
We introduce the following $\CC[f]$-modules $M_0$ and $M_1$:
$$
M_0=\frac{\{\omega\in \Omega^1(X), \exists \eta \in \Omega^1(X), d\omega = \eta\wedge df\}}{d\Omega^0(X)+ \Omega^0(X)df} 
$$
$$
M_1=\frac{\Omega^1(X)}{\{\omega\in \Omega^1(X), \exists \eta \in \Omega^1(X), d\omega = \eta\wedge df\}} 
$$
Note that we have the exact sequence of $\CC[f]$-modules:
$$
0\longrightarrow M_0 \longrightarrow H^1(f) \longrightarrow M_1 \longrightarrow 0
$$
Since localisation is an exact functor and $\CC[f]_{(P(f))}$ is a noetherian module for any $P\not=0$, it suffices to
prove that both $M_0$ and $M_1$ become finite modules after
a suitable localisation. The module $M_0$ is by definition the first group of standard relative cohomology of
$f$ (see \cite{}). By a theorem of Hartshorne (see \cite{Ha}), there exists a non-zero polynomial $P$ of $\CC[t]$ such that $(M_0)_{(P(f))}$
is a $\CC[f]_{(P(f))}$-module of finite type. So there only remains to prove that $M_1$ becomes a finite module after
a suitable localisation, and this is what we will do in the following lemmas.

\begin{lem} \label{finite2}
Let $X$ be an irreducible affine surface, $S$ its singular set and $I$ the defining ideal of $S$ in $\OX$. Let $f:X\rightarrow \CC$
be a dominant map. Then there exists
a non-zero polynomial $P$ of $\CC[t]$ such that, for any $n\geq 0$, the quotient $(\Omega^1(X)/I^n)_{(P(f))}$ is a finite
$\CC[f]_{(P(f))}$-module.
\end{lem}
\dem Since $S$ has dimension $\leq 1$, there exists a non-empty Zariski open set $U$ in $\CC$ such that either $f^{-1}(U)\cap S$
is empty or the restriction $f:f^{-1}(U)\cap S\longrightarrow U$ is a finite morphism. Let $P$ be a non-zero polynomial whose roots
form the set $\CC-U$. In the first case, $(\OS)_{(P(f))}$ is equal to zero. In the second case, the ring $(\OS)_{(P(f))}$ is a finite
$\CC[f]_{(P(f))}$-module. Since $I$ is radical, $\OS$ coincides with $\OX/I$ and $(\OX/I)_{(P(f))}$ is a finite
$\CC[f]_{(P(f))}$-module. It is then easy to prove that $(\OX/I^n)_{(P(f))}$ is a finite
$\CC[f]_{(P(f))}$-module, by an induction on $n$ and by using the following exact sequence:
$$
0\longrightarrow (I^n/I^{n+1})_{(P(f))} \longrightarrow (\OX/I^{n+1})_{(P(f))}\longrightarrow (\OX/I^n)_{(P(f))}\longrightarrow 0
$$
Here the only thing to note is that $I^n/I^{n+1}$ is a finite $\OX/I$-module for any $n$. Since $\Omega^1(X)$ is a finite $\OX$-module,
$\Omega^1(X)/I^n$ is a finite $\OX/I^n$-module. Therefore $(\Omega^1(X)/I^n)_{(P(f))}$ is a finite $\CC[f]_{(P(f))}$-module.
\qed

\begin{lem} \label{finite1}
Let $X$ be an irreducible affine surface, $S$ its singular set and $I$ the defining ideal of $S$ in $\OX$. Let $f:X\rightarrow \CC$
be a dominant map. Then there exists a non-zero polynomial $P$ of $\CC[t]$ and an integer $N$ such that:
$$
I^N \Omega^2(X)_{(P(f))} \subseteq \Omega^1(X)_{(P(f))} \wedge df
$$
\end{lem}
\dem By the generic smoothness theorem (see \cite{Jou}), there exists a non-empty Zariski open set $U$ of $\CC$ such that the restriction
$f: f^{-1}(U) \cap (X-S)\longrightarrow U$ is non-singular. Let $P$ be a non-zero polynomial whose roots
form the set $\CC-U$, and denote by $X'$ the surface $f^{-1}(U)\cap X$. We can then identify $\Omega^i(X)_{(P(f))}$
with $\Omega^i(X')$ for any $i$. We are going to prove there exists an integer $N$ such that:
$$
I^N \Omega^2(X') \subseteq \Omega^1(X') \wedge df
$$
Let $f_1,...,f_r$ be a set of generators of $I$. Let $\Omega$ be a regular 2-form on $X'$. Since $f_i$ belongs to $I$, the surface
$X'-V(f_i)$ is smooth and the restriction
$f: X'-V(f_i)\longrightarrow U$ is non-singular. By the De Rham lemma (see \cite{}), there exists a regular 1-form $\eta_i$ on
$X'-V(f_i)$ such that $\Omega= \eta_i \wedge df$ on $X'-V(f_i)$. Write $\eta_i$ as $\theta_i/f_i ^{n_i}$, where $\theta_i$ is regular
on $X'$. Then there exists an integer $m_i$ such that:
$$
f_i ^{m_i}(f_i ^{n_i}\Omega - \theta_i \wedge df)=0
$$
on $X'$, and so $f_i ^{m_i + n_i}\Omega$ belongs to $\Omega^1(X') \wedge df$. We set $N_{\Omega}=r\sup\{m_i +n_i\}$. Every element
$g$ of $I^{N_{\Omega}}$ can be written as a linear combination of the form:
$$
g=\sum{i_1 + ...+ i_r =N_{\Omega}} a_{i_1,...,i_r} f_1 ^{i_1} ...f_r ^{i_r}
$$
Since $i_1 + ...+ i_r =N_{\Omega}$, at least one of the indices $i_k$ is no less than $m_k + n_k$. So for any multi-index $(i_1,...,i_r)$,
$a_{i_1,...,i_r} f_1 ^{i_1} ...f_r ^{i_r} \Omega$ belongs to $\Omega^1(X') \wedge df$. Therefore $g\Omega$ belongs to $\Omega^1(X') \wedge df$
for any $g$ in $I^{N_{\Omega}}$. Now let $\Omega_1,...,\Omega_s$ be a set of generators of $\Omega^2(X')$ as an ${\cal{O}}_{X'}$-module. If
$N\geq N_{\Omega_i}$ for any $i=1,..,s$, then we have obviously:
$$
I^N \Omega^2(X') \subseteq \Omega^1(X') \wedge df
$$
\qed

\begin{lem}
Let $f:X\rightarrow \CC$ be a dominant map, where $X$ is an irreducible affine surface. Then there exists a non-zero polynomial $P$
of $\CC[t]$ such that $(M_1)_{(P(f))}$ is a finite $\CC[f]_{(P(f))}$-module.
\end{lem}
\dem We keep the notations of lemma \ref{finite1}. For any element $\omega$ of $I^{N+1}\Omega^1(X')$, the 2-form $\Omega=d\omega$
belongs to $I^N\Omega^2(X')$, hence to $\Omega^1(X')\wedge df$ by lemma \ref{finite1}. Therefore we have the inclusion:
$$
I^{N+1}\Omega^1(X')\subseteq \{\omega\in \Omega^1(X'), \exists \eta \in \Omega^1(X'), d\omega = \eta\wedge df\}
$$
By lemma \ref{finite2}, there exists a non-zero polynomial $Q$ such that $(\Omega^1(X)/I^n)_{(Q(f))}$ is a finite $\CC[f]_{(Q(f))}$-module.
The previous inclusion then induces the following surjective morphism of $\CC[f]_{(PQ(f))}$-modules:
$$
L: (\Omega^1(X)/I^n)_{(PQ(f))}\longrightarrow (M_1)_{(PQ(f))}
$$
Since $(\Omega^1(X)/I^n)_{(PQ(f))}$ is finite over $\CC[f]_{(PQ(f))}$, the result follows.
\qed

\subsection{Proof of proposition \ref{fini}}

In this subsection, we are going to prove more than proposition \ref{fini}. More precisely we are going to relate the rank of $H^1(f)$ (which is finite
by proposition \ref{finite3}) to the dimension of the $H^1(f^{-1}(y))$.

\begin{prop} \label{finite4}
Let $f:X\longrightarrow \CC$ be a dominant morphism, where $X$ is an irreducible affine surface. Let
$h_f$ be the rank of the module $H^1(f)$. Then for generic $y$ in $\CC$, the dimension of $H^1(f^{-1}(y))$
is equal to $h_f$.
\end{prop}

\begin{lem} \label{fini2}
Let $X$ be a complex affine surface and $f:X\longrightarrow \CC$ a dominant morphism.
If $(f-y)$ is a radical ideal in ${\cal{O}}_X$, then $H^1(f)/(f-y)\simeq H^1(f^{-1}(y))$. In particular,
this holds for generic $y$.
\end{lem}
\dem By definition, we have a first isomorphism:
\begin{equation*}
\begin{split}
H^1(f)/(f-y) &
\simeq \frac{\Omega^1(X)}{d\Omega^0(X)+\Omega^0(X)df+
(f-y)\Omega^1(X)} \\
&\simeq \frac{\Omega^1(X) /\Omega^0(X)df + (f-y)\Omega^1(X)}{
d\Omega^0(X)+ \Omega^0(X)df + (f-y)\Omega^1(X)/\Omega^0(X)df+(f-y)\Omega^1(X)}
\end{split}
\end{equation*}
Since $(f-y)$ is a radical ideal in ${\cal{O}}_X$, the restriction morphism induces an isomorphism:
$$
\Omega^1(X)/\Omega^0(X)df + (f-y)\Omega^1(X)\simeq \Omega^1(f^{-1}(y))
$$
From that we deduce $\displaystyle H^1(f)/(f-y) \simeq \Omega^1(f^{-1}(y))/d\Omega^0(f^{-1}(y))= H^1(f^{-1}(y))$.
\qed
{\it Proof of proposition \ref{finite4}}: Let $f:X\longrightarrow \CC$ be a dominant morphism, where $X$ is an
irreducible affine surface. Let $h_f$ be the rank of the module $H^1(f)$. By proposition \ref{finite3}, there
exists a non-zero polynomial $P$ of $\CC[t]$ such that $H^1(f)_{(P(f))}$ is a finite $\CC[f]_{(P(f))}$-module.
Up to refining the localisation, we may even assume that $H^1(f)_{(P(f))}$ is a finite free $\CC[f]_{(P(f))}$-module
of rank $h_f$. For any $y$ such that $(f-y)$ is a radical ideal in $\OX$ and $P(y)\not=0$, we have by lemma \ref{fini2}:
$$
H^1(f)_{(P(f))}/(f-y)\simeq H^1(f)/(f-y)\simeq H^1(f^{-1}(y))
$$
Since $H^1(f)_{(P(f))}$ is finite free of rank $h_f$, $H^1(f)_{(P(f))}/(f-y)$ has dimension $h_f$ and the result follows.

\section{The property ${\cal{P}}$}

In this subsection, we are going to prove the inequality given in theorem \ref{second}
by using a special property of the relative cohomology group $H^1(f)$. This
property will enable us to control the dimension of $H^1(f^{-1}(t))$ by means of the rank of $H^1(f)$.

\begin{df}
A $\CC[f]$-module $M$ satisfies the property ${\cal{P}}(y)$ if for any integer
$r$ and any element $\omega$ of $M$, we have: $(f-y)^r \omega=0 \Longrightarrow
\omega \in (f-y)M$.
\end{df}

\begin{lem} \label{alg}
Let $M$ be a $\CC[f]$-module satisfying ${\cal{P}}(y)$. Then $dim \; M/(f-y) \leq rk \; M$.
\end{lem}
\dem Let $e_1,..,e_s$ be some elements of $M$
whose classes in $M/(f-y)$ are free. In order to establish the lemma, we prove by contradiction
that $e_1,..,e_s$ are free in $M$. Assume there exist some polynomials $P_1(f),..,P_s(f)$ not all zero
such that $P_1(f)e_1+...+P_s(f)e_s=0$ in $M$. Let $m$ be the minimum of the orders of the $P_i$ at
$y$. Every $P_i(f)$ can be written as $P_i(f)=(f-y)^m T_i(f)$
where at least one of the $T_i(y)$ is nonzero. So we get:
$$
(f-y)^m\left \{ T_1(f)e_1+...+T_s(f)e_s\right \}=0
$$
By the property ${\cal{P}}(y)$, this implies:
$$
T_1(f)e_1+...+T_s(f)e_s\equiv T_1(y)e_1+...+T_s(y)e_s \equiv 0\; [(f-y)]
$$
Since the $e_i$ are free modulo $(f-y)$, every $T_i(y)$ is zero, hence a contradiction.
\qed
Our purpose in this subsection is to prove: 
 
\begin{prop} \label{PP}
Let $X$ be a complex irreducible affine surface, and $f:X\longrightarrow \CC$ a dominant morphism.
Assume that $X$ is locally a complete intersection. If $f^{-1}(y)\cap Sing(f)$ is finite, then
$H^1(f)$ satisfies the property ${\cal{P}}(y)$.
\end{prop}
Since $X$ is locally a complete intersection, the finiteness of $f^{-1}(y)\cap Sing(f)$ implies
that $(f-y)$ is a radical ideal in ${\cal{O}}_X$. By lemma \ref{fini2}, we have $H^1(f)/(f-y)\simeq H^1(f^{-1}(y))$.
So theorem \ref{second} will follow from lemma \ref{alg} and proposition \ref{PP}. We begin
with a few lemmas.

\begin{lem} \label{PP2}
Let $X$ be a complex affine surface that is locally a complete intersection.
Let $\omega$ be a regular
1-form on $X$ and $A$ a regular function on $X$ such that $(f-y)\omega=Adf$. If $f^{-1}(y)\cap Sing(f)$ is finite, there exists
a regular function $B$ on $X$ such that $\omega=Bdf$.
\end{lem}
\dem Let $\omega$ be a regular 1-form on $X$ and $A$ a regular function on $X$ such that
$(f-y)\omega=Adf$. Then $A$ vanishes on the set $f^{-1}(y)-Sing(f)$. Since $f^{-1}(y)$
is equidimensionnal of dimension 1 and $f^{-1}(y)\cap Sing(f)$ is finite, $A$ vanishes
on $f^{-1}(y)$. Since $f^{-1}(y)\cap Sing(f)$ is finite and $X$ is locally a complete
intersection, $f^{-1}(y)$ defines locally a complete intersection. Hence it is a complete
intersection on $X$, and $(f-y)$ divides $A$. If $A=(f-y)B$, then $(f-y)(\omega -Bdf)=0$.
Since $X$ is locally a complete intersection, the module $\Omega^1(X)$ is torsion-free (see
\cite{Gr}) and $\omega =Bdf$.
\qed

\begin{lem}
Let $X$ be a complex irreducible affine surface and $f:X\rightarrow \CC$ a dominant morphism. Let $C_1,..,C_r$ 
be the connected components of $f^{-1}(t)$ and $n$ an integer $\geq 0$. Then there exist some regular functions
$S_{i,n}$ on $X$ such that $S_{i,n}=1$ on $C_i$, $S_{i,n}=0$ on $C_j$ for $j\not=i$ and
$dS_{i,n}$ belongs to $(f-t)^{n+1}\Omega^1(X)$.
\end{lem}
\dem For simplicity, assume that $t=0$. There exists a regular function $T_i$ on $X$ such that
$T_i=1$ on $C_i$ and $T_i=0$ on $C_j$ for $j\not=i$. Then $T_i(1-T_i)$ vanishes
on $f^{-1}(0)$ and by Hilbert's Nullstellensatz, there exists an integer $m$ such that
$T_i ^m(1-T_i)^m$ belongs to $f^{n+1}{\cal{O}}_X$. We set:
$$
P_i(x)=\int_{0} ^{x} t^m(1-t)^m dt \quad \mbox{and} \quad R_{i,n} = P_i(T_i)
$$
By construction the 1-form $dR_{i,n}=T_i ^m(1-T_i)^m dT_i$ is divisible by $f^{n+1}$. Since
$P_i(0)=0$ and $T_i$ vanishes on $C_j$ for $j\not=i$, $R_{i,n}$ vanishes on $C_j$ if $j\not=0$. Since $P_i(1)\not=0$,
$R_{i,n}=P_i(1)\not=0$ on $C_i$. Then choose $S_{i,n}=R_{i,n}/P_i(1)$.
\qed

\begin{lem} \label{PP1}
Let $X$ be a complex irreducible affine surface and $f:X\rightarrow \CC$ a dominant morphism.
Let $R$ be a regular function on $X$ such that $dR =Adf + (f-t)\eta$, where $A,\eta$ are regular on $X$.
Then $R$ is locally constant on $f^{-1}(t)$. 
\end{lem}
\dem Since $dR =Adf + (f-t)\eta$, the restriction of $dR$ to $f^{-1}(t)$ is zero. So $R$
is singular at any smooth point of $f^{-1}(t)$, and $R$ is constant on every
connected component of the smooth part of $f^{-1}(t)$. By continuity and density,
$R$ is constant on every connected component of $f^{-1}(t)$, hence it is locally
constant on $f^{-1}(t)$.
\qed
{\it{Proof of proposition \ref{PP}:}} Let $X$ be a complex irreducible affine surface that is locally a complete intersection.
Let $f:X\rightarrow \CC$ be a dominant morphism and assume that $f^{-1}(t)\cap Sing (f)$ is finite. We may assume that $t=0$.
Let us prove by induction on $n\geq 0$ that, if $f^n\omega=0$ in $H^1(f)$, then $\omega$ belongs to $(f)H^1(f)$. This is trivial
for $n=0$. Assume that the assertion holds to the order $n$. Let $\omega$
be a regular 1-form on $X$ such that $f^{n+1}\omega=0$ in $H^1(f)$. Then
there exist some regular functions $R,A$ such that $f^{n+1}\omega=dR +Adf$
on $\Omega^1(X)$. By lemma \ref{PP1}, $R$ is locally constant on $f^{-1}(0)$.
Let $C_1,..,C_r$ be the connected components of $f^{-1}(0)$. If $R$ takes the value $\lambda_i$ on $C_i$, then the function:
$$
R'= R - \sum_i \lambda_i S_{i,n+1}
$$
vanishes on $f^{-1}(0)$. By construction, there exists a regular 1-form $\eta$ such that:
$$
f^{n+1}\omega=dR' +Adf + f^{n+2}\eta
$$
Since $f^{-1}(0)\cap Sing(f)$ is finite and $X$ is locally a complete intersection, $(f)$ is a radical ideal and $R'$ is divisible by $f$.
If $R'=fS$ with $S$ regular, we obtain:
$$
f\left(f^{n}\omega- dS - f^{n+1}\eta \right)=(A+S)df
$$
By lemma \ref{PP2}, there exists a regular function $B$ such that:
$$
f^{n}(\omega -f \eta) = dS +Bdf
$$
By induction $(\omega-f\eta)$ belongs to $(f)H^1(f)$, as well as $\omega$, and we are done.
\qed

\section{An example}

We end this paper with an example of a surface that is not locally a complete intersection (for more details,
see \cite{Di}).
For that surface there exists a map for which the conclusion of theorem
\ref{second}
fails. Let $(u,v,w_1, w_2)$ be a system of coordinates in $\CC^4$, and consider the affine
set $X$ of $\CC^4$ defined by the equations:
$$
u^2 w_1  -v^2=0, \quad u^3 w_2 -v^3=0, \quad w_1 ^3 - w_2 ^2=0
$$
Note that $X$ can be reinterpreted as:
$$
X=Spec(\CC[x,xy,y^2,y^3])
$$
So $X$ is an irreducible surface. Moreover 0 is the only
singular point of $X$, but $X$ is not locally a complete intersection.
Indeed if it were so, then $X$ would be a normal surface because
it
is non-singular in codimension 1. Consider the function
$h= w_2/w_1=v/u$ on $X$. It
is well-defined and regular outside the origin, hence $h$ is regular
because $X$ is normal. Moreover we have the following relations:
$$
v=hu, \quad \quad w_1= h^2, \quad \quad w_2=h^3
$$
So every regular function on $X$ can be expressed as a
polynomial in $(u,h)$, and $X$ is isomorphic to $\CC^2$.
But this is impossible because $X$ is singular at the origin.
Consider now the map
$f: X\rightarrow \CC$ defined by:
$$
f(u,v,w_1,w_2)=u
$$
For
$y\not=0$, the fibre $f^{-1}(y)$ is isomorphic to a line, hence
$H^1(f^{-1}(y))=0$. The fibre $f^{-1}(0)$ is isomorphic
to a cusp, hence contractible, and $f^{-1}(0)\cap Sing(f)$ is reduced to
the origin. Moreover its Milnor number coincides
with its local Betti number and is equal to 2. With the notations
of the previous sections, $h_f=0$ and $\dim \; H^1(f^{-1}(0))=2$,
so that $\dim \; H^1(f^{-1}(0))\not\leq h_f$.

\end{document}